\documentclass[a4paper,11pt]{amsart}
\usepackage{graphicx}
\usepackage{mathptmx}
\usepackage{amsmath}
\usepackage{amssymb}
\usepackage{enumitem}
\usepackage{xcolor}

\newmuskip\pFqmuskip

\newcommand*\pFq[6][8]{%
  \begingroup 
  \pFqmuskip=#1mu\relax
  \mathcode`=\string"8000
  \begingroup\lccode`\~=`\,
  \lowercase{\endgroup\let~}\pFqcomma
  F^{#2}_{#3}{\left(\genfrac..{0pt}{}{#4}{#5}\bigg|#6\right)}%
  \endgroup
}
\newcommand{\pFqcomma}{\mskip\pFqmuskip}

\newtheorem{theorem}{Theorem}

\newtheorem{corollary}[theorem]{Corollary}

\begin{document}

\title[Representations of degenerate poly-Bernoulli polynomials]{Representations of degenerate poly-Bernoulli polynomials}

\author{${}^{1}$Taekyun  Kim}
\address{$^{1}$Department of Mathematics, Kwangwoon University, Seoul 139-701, Republic of Korea}
\email{tkkim@kw.ac.kr}

\author{${}^{2}$DAE SAN KIM}
\address{$^{2}$Department of Mathematics, Sogang University, Seoul 121-742, Republic of Korea}
\email{dskim@sogang.ac.kr}

\author{${}^{3}$Jongkyum Kwon}
\address{${3}$Department of Mathematics Education and ERI, Gyeongsang National University Jinju 660-701, Repulbic of Korea}
\email{mathkjk26@gnu.ac.kr}

\author{${}^{4}$Hyunseok Lee}
\address{$^{4}$Department of Mathematics, Kwangwoon University, Seoul 139-701, Republic of Korea}
\email{luciasconstant@kw.ac.kr}


\begin{abstract}
As is well-known, poly-Bernoulli polynomials are defined in terms of polylogarithm functions. Recently, as degenerate versions of such functions and polynomials, degenerate polylogarithm functions were introduced and degenerate poly-Bernoulli polynomials were defined by means of the degenerate polylogarithm functions, and some of their properties were investigated.
The aim of this paper is to further study some properties of the degenerate poly-Bernoulli polynomials by using three formulas coming from the recently developed `$\lambda$-umbral calculus.' In more detail, among other things, we represent the degenerate poly-Bernoulli polynomials by higher-order degenerate Bernoulli polynomials and by higher-order degenerate derangement polynomials.
\end{abstract}

\maketitle

\section{Introduction}

Carlitz is the first one who initiated the study of degenerate versions of some special numbers and polynomials, namely the degenerate Bernoulli and Euler polynomials and numbers (see [1]).
In recent years, studying degenerate versions of some special numbers and polynomials regained interests of some mathematicians with their interests not only in combinatorial and arithmetic properties but also in applications to differential equations, identities of symmetry and probability theory (see [5,6,9,10,13,15,16] and the references therein). It is noteworthy that studying degenerate versions is not only limited to polynomials but also can be extended to transcendental functions like gamma functions (see [8]). \par
The Rota's theory of umbral calculus is based on linear functionals and differential operators (see [2-4,17-21]). The Sheffer sequences occupy the central position in the theory and are characterized by the generating functions where the usual exponential function enters. The motivation for the paper [6] starts from the question that what if the usual exponential function is replaced by the degenerate exponential functions (see \eqref{2}). As it turns out, it corresponds to replacing the linear functional by the family of $\lambda$-linear functionals (see \eqref{12}) and the differential operator by the family of  $\lambda$-differential operators (see \eqref{15}). Indeed, these replacements lead us to define $\lambda$-Sheffer polynomials and degenerate Sheffer polynomials (see \eqref{17}). \par
As is well-known, poly-Bernoulli polynomials are defined in terms of polylogarithm functions. Recently, as degenerate versions of such functions and polynomials, degenerate polylogarithm functions were introduced and degenerate poly-Bernoulli polynomials were defined by means of the degenerate polylogarithm functions, and some properties of the degenerate poly-Bernoulli polynomials were investigated (see [13]). \par
The aim of this paper is to further study the degenerate poly-Bernoulli polynomials, which is a $\lambda$-Sheffer sequence and hence a degenerate Sheffer sequence, by using the above-mentioned $\lambda$-linear functionals and $\lambda$-differential operators.
In more detail, these polynomials are investigated by three different tools, namely a formula about representing a $\lambda$-Sheffer sequence by another (see \eqref{20}), a formula coming from the generating functions of  $\lambda$-Sheffer sequences (see Theorem 1) and a formula arising from the definitions for $\lambda$-Sheffer sequences (see Theorems 6,7). Then, among other things, we represent the degenerate poly-Bernoulli polynomials by higher-order degenerate Bernoulli polynomials and by higher-order degenerate derangement polynomials. The rest of this section is devoted to recalling the necessary facts that are needed throughout the paper, which includes the `$\lambda$-umbral calculus'.

\vspace{0.2cm}

For $k\in\mathbb{Z}$, and $0 \ne \lambda\in\mathbb{R}$, the degenerate polylogarithm functions are defined by
\begin{equation}
\mathrm{Li}_{k,\lambda}(x)=\sum_{n=1}^{\infty}\frac{(-\lambda)^{n-1}(1)_{n,1/\lambda}}{(n-1)!n^{k}}x^{n},\quad(\mathrm{see}\ [5]),\label{1}
\end{equation}
 where $(x)_{0,\lambda}=1,\ (x)_{n,\lambda}=x(x-\lambda)\cdots(x-(n-1)\lambda)$, $(n\ge 1)$. \par
For any $\lambda\in\mathbb{R}$, the degenerate exponential functions are given by
 \begin{equation}
 e_{\lambda}^{x}(t)=\sum_{l=0}^{n}\frac{(x)_{n,\lambda}}{n!}t^{n},\quad e_{\lambda}(t)=e_{\lambda}^{1}(t)=\sum_{n=0}^{\infty}\frac{(1)_{n,\lambda}}{n!}t^{n},\quad(\mathrm{see}\ [8]).\label{2}
\end{equation}
The compositional inverse $\log_{\lambda}(t)$ of $e_{\lambda}(t)$ is given by
\begin{equation}
\log_{\lambda}(t)=\sum_{n=1}^{\infty}\frac{\lambda^{n-1}}{n!}(1)_{n,1/\lambda}(t-1)^{n},\quad(\mathrm{see}\ [5]).\label{3}
\end{equation}
From \eqref{1} and \eqref{3}, we have
\begin{equation}
\mathrm{Li}_{1,\lambda}(x)=-\log_{\lambda}(1-x),\quad\mathrm{and}\quad\lim_{\lambda\rightarrow 0}\mathrm{Li}_{k,\lambda}(x)=\mathrm{Li}_{k}(x),\label{4}
\end{equation}
where $\mathrm{Li}_{k}(x)$ are the polylogarithm functions defined by
\begin{displaymath}
\mathrm{Li}_{k}(x)=\sum_{n=1}^{\infty}\frac{x^{n}}{n^{k}},\quad(\mathrm{see}\ [5,7]).
\end{displaymath}
In [5], the degenerate poly-Bernoulli polynomials are defined by Kim-Kim as
\begin{equation}
\frac{\mathrm{Li}_{k,\lambda}(1-e_{\lambda}(-t))}{e_{\lambda}(t)-1}e_{\lambda}^{x}(t)=\sum_{n=0}^{\infty}B_{n,\lambda}^{(k)}(x)\frac{t^{n}}{n!}.\label{5}
\end{equation}
When $x=0$, $B_{n,\lambda}^{(k)}= B_{n,\lambda}^{(k)}(0)$ are called the degenerate poly-Bernoulli numbers.\par
It is well known that Carlitz's degenerate Bernoulli polynomials of order $r$ are defined by
\begin{equation}
\bigg(\frac{t}{e_{\lambda}(t)-1}\bigg)^{r}e_{\lambda}^{x}(t)=\sum_{n=0}^{\infty}\beta_{n,\lambda}^{(r)}(x)\frac{t^{n}}{n!},\quad(\mathrm{see}\ [1]).\label{6}
\end{equation}
For $r=1$, $\beta_{n,\lambda}^{(r)}(x)=\beta_{n,\lambda}(x)$ are called the degenerate Bernoulli polynomials. \par
From \eqref{5}, we note that
\begin{equation}
\sum_{n=0}^{\infty}B_{n,\lambda}^{(1)}(x)\frac{t^{n}}{n!}=\frac{\mathrm{Li}_{1,\lambda}(1-e_{\lambda}(-t))}{e_{\lambda}(t)-1}   =\frac{t}{e_{\lambda}(t)-1}e_{\lambda}^{x}(t)=\sum_{n=0}^{\infty}\beta_{n,\lambda}(x)\frac{t^{n}}{n!}. \label{7}
\end{equation}
By \eqref{7}, we get $B_{n,\lambda}^{(1)}(x)=\beta_{n,\lambda}(x)$, $(n\ge 0)$, (see [7]). \par
The degenerate Stirling numbers of the second kind appear as the coefficients in the expansion
\begin{equation}
(x)_{n,\lambda}=\sum_{l=0}^{n}S_{2,\lambda}(n,l)(x)_{l},\quad(n\ge 0),\quad(\mathrm{see}\ [5]).\label{8}
\end{equation}
As the inversion formula of \eqref{8}, the degenerate Stirling numbers of first kind appear as the coefficients in the expansion
\begin{equation}
(x)_{n}=\sum_{l=0}^{n}S_{1,\lambda}(n,l)(x)_{l,\lambda},\quad(n\ge 0),\quad(\mathrm{see}\ [5]).\label{9}
\end{equation}
Thus, by \eqref{8} and \eqref{9}, we have
\begin{equation}
\frac{1}{k!}\big(e_{\lambda}(t)-1\big)^{k}=\sum_{n=k}^{\infty}S_{2,\lambda}(n,k)\frac{t^{n}}{n!},\label{10}
\end{equation}
and
\begin{displaymath}
\frac{1}{k!}\big(\log_{\lambda}(1+t)\big)^{k}=\sum_{n=k}^{\infty}S_{1,\lambda}(n,k)\frac{t^{n}}{n!},\quad(k\ge 0),\quad(\mathrm{see}\ [5,9,11,16]).
\end{displaymath}
In view of \eqref{6}, the degenerate derangement polynomials of order $r(\in\mathbb{N})$ are defined by
\begin{equation}
\frac{1}{(1-t)^{r}}e_{\lambda}^{-1}(t) e_{\lambda}^{x}(t)=\sum_{n=0}^{\infty}d_{n,\lambda}^{(r)}(x)\frac{t^{n}}{n!},\quad(\mathrm{see}\ [10]).\label{11}
\end{equation}
When $r=1$, $d_{n,\lambda}(x)=d_{\lambda}^{(1)}(x)$ are called the degenerate derangement polynomials. \par
Note that $\displaystyle \lim_{\lambda\rightarrow 0}d_{n,\lambda}(x)=d_{n}(x)\displaystyle$, where $d_{n}(x)$ are the derangement polynomials and $d_{n}=d_{n}(0)$ the derangement numbers (see [12,13,14,15]). \par
We remark that the umbral calculus has long been studied by many people (see [2,3,4,6,18-20]). For the rest of this section, we will recall the necessary facts on the $\lambda$-linear functionals, $\lambda$-differential operators and $\lambda$-Sheffer sequences, and so on. The details on these can be found in the recent paper [6]. \par
Let $\mathbb{C}$ be the field of complex numbers,
\begin{displaymath}
\mathcal{F}=\bigg\{f(t)=\sum_{k=0}^{\infty}a_{k}\frac{t^{k}}{k!}\ \bigg|\ a_{k}\in\mathbb{C}\bigg\},
\end{displaymath}
and let
\begin{displaymath}
\mathbb{P}=\mathbb{C}[x]=\bigg\{\sum_{i=0}^{\infty}a_{i}x^{i}\ \bigg|\ a_{i}\in\mathbb{C}\textrm{ with $a_{i}=0$ for all but finite number of $i$}\bigg\}.
\end{displaymath}
For $f(t)\in\mathcal{F}$ with $\displaystyle f(t)=\sum_{k=0}^{\infty}a_{k}\frac{t^{k}}{k!}\displaystyle$, and $\lambda\in\mathbb{R}$, the $\lambda$-linear functional $\langle f(t)|\cdot\rangle_{\lambda}$ on $\mathbb{P}$ is defined by
\begin{equation}
\langle f(t)|(x)_{n,\lambda}\rangle_{\lambda}=a_{n},\quad(n\ge 0),\quad(\mathrm{see}\ [6]). \label{12}
\end{equation}
By \eqref{12}, we get
\begin{equation}
\langle t^{k}|(x)_{n,\lambda}\rangle=n!\delta_{n,k},\quad(n,k\ge 0),\quad(\mathrm{see}\ [6]),\label{13}
\end{equation}
where $\delta_{n,k}$ is the Kronecker's symbol. \par
The $\lambda$-differential operators on $\mathbb{P}$ are defined by
\begin{equation}
(t^{k})_{\lambda}(x)_{n,\lambda}=\left\{\begin{array}{ccc}
(n)_{k}(x)_{n-k,\lambda}, & \textrm{if $0 \le k\le n$},\\
0, & \textrm{if $k>n$.}
\end{array}\right.\label{15}
\end{equation}
For $\displaystyle f(t)=\sum_{k=0}^{\infty}a_{k}\frac{t^{k}}{k!}\in\mathcal{F}\displaystyle$, and by \eqref{15}, we get
\begin{align}
\big(f(t)\big)_{\lambda}(x)_{n,\lambda}\ &=\ \sum_{k=0}^{n}\binom{n}{k}a_{k}(x)_{n-k,\lambda},\quad(n\ge 0), \label{16}\\
\big(e_{\lambda}^{y}(t)\big)_{\lambda}(x)_{n,\lambda}\ &=\ (x+y)_{n,\lambda},\quad(n\ge 0),\quad(\mathrm{see}\ [6]).\nonumber
\end{align} \par
Let $f(t)$ be a delta series and let $g(t)$ be an invertible series. Then there exists a unique sequence $S_{n,\lambda}(x)\ (\deg S_{n,\lambda}(x)=n)$ of polynomials satisfying the orthogonality conditions
\begin{equation}
\Big\langle g(t)\big(f(t)\big)^{k}\ \Big|\ S_{n,\lambda}(x)\Big\rangle_{\lambda}=n!\delta_{n,k},\quad(n,k\ge 0),\quad(\mathrm{see}\ [6]). \label{17}
\end{equation}
Here $S_{n,\lambda}(x)$ is called the $\lambda$-Sheffer sequence for $\big(g(t), f(t)\big)$, which is denoted by $S_{n,\lambda}(x) \sim \big(g(t),f(t)\big)_{\lambda}$.
The sequence $S_{n,\lambda}(x)$ is the $\lambda$-Sheffer sequence for $(g(t),f(t))$ if and only if
\begin{equation}
\frac{1}{g(\overline{f}(t))}e_{\lambda}^{y}\big(\overline{f}(t)\big)=\sum_{n=0}^{\infty}S_{n,\lambda}(y)\frac{t^{n}}{n!},\quad(\mathrm{see}\ [6]),\label{18}
\end{equation}
for all $y\in\mathbb{C}$, where $\overline{f}(t)$ is the compositional inverse of $f(t)$ such that $f(\overline{f}(t))=\overline{f}(f(t))=t$. \par
Let $S_{n,\lambda}(x)\sim\big(g(t),f(t)\big)_{\lambda}$. Then, from Theorem 16 of [6], we recall that
\begin{equation}
 (f(t))_{\lambda}S_{n,\lambda}(x)=nS_{n-1,\lambda}(x), (n\ge 1). \label{19}
\end{equation}
For $S_{n,\lambda}(x)\sim (g(t),f(t))_{\lambda}$, $r_{n,\lambda}(x)\sim (h(t),l(t))$, we have
\begin{equation}
S_{n,\lambda}(x)=\sum_{k=0}^{m}C_{n,k}r_{k,\lambda}(x),\quad(n\ge 0),\quad(\mathrm{see}\ [6]),\label{20}
\end{equation}
where
\begin{displaymath}
C_{n,k}=\frac{1}{k!}\bigg\langle\frac{h(\overline{f}(t))}{g(\overline{f}(t))}\big(l(\overline{f}(t))\big)^{k}\ \bigg|\ (x)_{n,\lambda}\bigg\rangle_{\lambda}.
\end{displaymath}
In this paper, we study the properties of degenerate poly-Bernoulli polynomial arising from degenerate polylogarithmic function and give some identities of those polynomials associated with special polynomials which are derived from the properties of $\lambda$-Sheffer sequences.

\section{Representations of degenerate poly-Bernoulli polynomials}
For $S_{n,\lambda}(x)\sim (g(t),f(t))_{\lambda},\ (n\ge 0)$, we have
\begin{align}
\Bigg\langle\frac{1}{g(\overline{f}(t))}e_{\lambda}^{x}\big(\overline{f}(t)\big)\ \Bigg|\ (x)_{n,\lambda}\Bigg\rangle_{\lambda}\ &=\ \sum_{k=0}^{\infty}S_{k,\lambda}(x)\frac{1}{k!}\big\langle t^{k}|(x)_{n,\lambda}\big\rangle_{\lambda} \label{21}\\
&=\ S_{n,\lambda}(x),\quad(n\ge 0).\nonumber
\end{align}
From \eqref{21}, we note that
\begin{align}
S_{n,\lambda}(x)\ &=\ \Bigg\langle\frac{1}{g(\overline{f}(t))}e_{\lambda}^{x}\big(\overline{f}(t)\big)\ \Bigg|\ (x)_{n,\lambda}\Bigg\rangle_{\lambda}\ =\ \sum_{j=0}^{\infty}\frac{1}{j!}\Bigg\langle\frac{1}{g(\overline{f}(t)}\big(\overline{f}(t)\big)^{j}\ \Bigg|\ (x)_{n,\lambda}\Bigg\rangle_{\lambda}(x)_{j,\lambda}\label{22}\\
&=\ \sum_{j=0}^{n}  \frac{1}{j!}\Bigg\langle\frac{1}{g(\overline{f}(t)}\big(\overline{f}(t)\big)^{j}\ \Bigg|\ (x)_{n,\lambda}\Bigg\rangle_{\lambda}(x)_{j,\lambda}.\nonumber
\end{align}
Therefore, by \eqref{22}, we obtain the following theorem.
\begin{theorem}
For $S_{n,\lambda}(x)\sim(g(t),f(t))_{\lambda}$, we have
\begin{displaymath}
S_{n,\lambda}(x)=\sum_{j=0}^{n} \frac{1}{j!}\Bigg\langle\frac{1}{g(\overline{f}(t)}\big(\overline{f}(t)\big)^{j}\ \Bigg|\ (x)_{n,\lambda}\Bigg\rangle_{\lambda}(x)_{j,\lambda}.
\end{displaymath}
\end{theorem}
From \eqref{5} and \eqref{18}, we have
\begin{equation}
B_{n,\lambda}^{(k)}(x)\sim \bigg(\frac{e_{\lambda}(t)-1}{\mathrm{Li}_{k,\lambda}(1-e_{\lambda}(-t))},t\bigg)_{\lambda}. \label{23}
\end{equation}
By Theorem 1, we get the following corollary.
\begin{corollary}
For $n\ge 0$, we have
\begin{displaymath}
B_{n,\lambda}^{(k)}(x)=\sum_{j=0}^{n}\binom{n}{j}B_{n-j,\lambda}^{(k)}(x)_{j,\lambda}.
\end{displaymath}
\end{corollary}
From \eqref{23}, we note that
\begin{equation}
\bigg(\frac{e_{\lambda}(t)-1}{\mathrm{Li}_{k,\lambda}(1-e_{\lambda}(-t))}\bigg)_{\lambda}B_{n,\lambda}^{(k)}(x)=(x)_{n,\lambda}\sim (1,t)_{\lambda}. \label{24}
\end{equation}
By \eqref{24}, and noting that $(x)_{n,\lambda}=\sum_{l=0}^{n}\sum_{j=0}^{l}S_{2,\lambda}(n,l)S_{1,\lambda}(l,j)(x)_{j,\lambda}$, we get
\begin{align}
B_{n,\lambda}^{(k)}(x)\ &=\ \bigg(\frac{\mathrm{Li}_{k,\lambda}(1-e_{\lambda}(-t))}{e_{\lambda}(t)-1}\bigg)_{\lambda}(x)_{n,\lambda}\label{25} \\
&=\ \sum_{l=0}^{n}\sum_{j=0}^{l}S_{2,\lambda}(n,l)S_{1,\lambda}(l,j)\bigg(\frac{\mathrm{Li}_{k,\lambda}(1-e_{\lambda}(-t))}{e_{\lambda}(t)-1}\bigg)_{\lambda}(x)_{j,\lambda}\nonumber \\
&=\ \sum_{l=0}^{n}\sum_{j=0}^{l}S_{2,\lambda}(n,l)S_{1,\lambda}(l,j)\sum_{m=0}^{j}B_{m,\lambda}^{(k)}\frac{1}{m!}(t^{m})_{\lambda}(x)_{j,\lambda} \nonumber \\
&=\ \sum_{l=0}^{n}\sum_{j=0}^{l}\sum_{m=0}^{j}\binom{j}{m}S_{2,\lambda}(n,l)S_{1,\lambda}(l,j)B_{m,\lambda}^{(k)}(x)_{j-m,\lambda} \nonumber\\
&=\ \sum_{l=0}^{n}\sum_{m=0}^{l}\sum_{j=m}^{l}\binom{j}{m}S_{2,\lambda}(n,l)S_{1,\lambda}(l,j)B_{m,\lambda}^{(k)}(x)_{j-m,\lambda} \nonumber \\
&=\ \sum_{l=0}^{n}\sum_{m=0}^{l}\sum_{j=0}^{l-m}\binom{j+m}{m}S_{2,\lambda}(n,l)S_{1,\lambda}(l,j+m)B_{m,\lambda}^{(k)}(x)_{j,\lambda}. \nonumber
\end{align}
Therefore, by \eqref{25}, we obtain the following theorem.
\begin{theorem}
For $n\ge 0$, we have
\begin{displaymath}
B_{n,\lambda}^{(k)}(x)= \sum_{l=0}^{n}\sum_{m=0}^{l}\sum_{j=0}^{l-m}\binom{j+m}{m}S_{2,\lambda}(n,l)S_{1,\lambda}(l,j+m)B_{m,\lambda}^{(k)}(x)_{j,\lambda}.
\end{displaymath}
\end{theorem}
Now, we observe that
\begin{align}
B_{n,\lambda}^{(k)}(y)\ &=\ \Bigg\langle \frac{\mathrm{Li}_{k,\lambda}(1-e_{\lambda}(-t))}{e_{\lambda}(t)-1}e_{\lambda}^{y}(t)\ \Bigg|\ (x)_{n,\lambda}\Bigg\rangle_{\lambda}\label{27}\\
&=\ \Bigg\langle \frac{\mathrm{Li}_{k,\lambda}(1-e_{\lambda}(-t))}{t}\ \Bigg|\ \bigg(\frac{t}{e_{\lambda}(t)-1}e_{\lambda}^{y}(t)\bigg)_{\lambda}(x)_{n,\lambda}\Bigg\rangle_{\lambda}\nonumber  \\
&=\ \sum_{l=0}^{n}\binom{n}{l}\beta_{l,\lambda}(y)\bigg\langle\frac{1}{t}\mathrm{Li}_{k,\lambda}\big(1-e_{\lambda}(-t)\big)\ \bigg|\ (x)_{n-l,\lambda}\bigg\rangle_{\lambda} \nonumber
\end{align}
\begin{align}
&=\ \sum_{l=0}^{n}\binom{n}{l}\beta_{l,\lambda}(y)\Bigg\langle\frac{1}{t}\sum_{m=1}^{\infty}\frac{(-\lambda)^{m-1}(1)_{m,1/\lambda}}{m^{k-1}}\frac{(-1)^{m}}{m!}\big(e_{\lambda}(-t)-1\big)^{m}\ \Bigg|\ (x)_{n-l,\lambda}\Bigg\rangle_{\lambda}\nonumber \\
&=\ \sum_{l=0}^{n}\binom{n}{l}\beta_{l,\lambda}^{(y)}\Bigg\langle \frac{1}{t}\sum_{j=1}^{\infty}\bigg(\sum_{m=1}^{j}\frac{(-\lambda)^{m-1}(1)_{m,1/\lambda}}{m^{k-1}} (-1)^{j-m}S_{2,\lambda}(j,m)\frac{t^{j}}{j!}\ \Bigg|\ (x)_{n-l,\lambda}\Bigg\rangle_{\lambda} \nonumber \\
&=\ \sum_{l=0}^{n}\binom{n}{l}\beta_{l,\lambda}(y)\sum_{j=0}^{\infty}\bigg(\sum_{m=1}^{j+1}\frac{(-\lambda)^{m-1}(1)_{m,1/\lambda}}{m^{k-1}}(-1)^{j+1-m}\frac{S_{2,\lambda}(j+1,m)}{(j+1)!}\bigg)\langle t^{j}|(x)_{n-l,\lambda}\rangle_{\lambda}\nonumber  \\
&=\ \sum_{l=0}^{n}\binom{n}{l}\beta_{l,\lambda}(y)\sum_{m=1}^{n-l+1}\frac{(-\lambda)^{m-1}(1)_{m,1/\lambda}}{m^{k-1}}\frac{(-1)^{n-l+1-m}}{(n-l+1)!}S_{2,\lambda}(n-l+1,m)(n-l)! \nonumber \\
&=\ \sum_{l=0}^{n}\sum_{m=1}^{n-l+1}\frac{(-1)^{n-l}\binom{n}{l}\lambda^{m-1}}{(n-l+1)m^{k-1}}(1)_{m,1/\lambda}S_{2,\lambda}(n-l+1,m)\beta_{l,\lambda}(y).\nonumber
\end{align}
Therefore, by \eqref{27}, we obtain the following theorem.
\begin{theorem}
For $n\ge 0$, we have
\begin{displaymath}
B_{n,\lambda}^{(k)}(x)=\sum_{l=0}^{n}\sum_{m=1}^{n-l+1}\frac{(-1)^{n-l}\binom{n}{l}\lambda^{m-1}}{(n-l+1)m^{k-1}}(1)_{m,1/\lambda}S_{2,\lambda}(n-l+1,m)\beta_{l,\lambda}(x).
\end{displaymath}
\end{theorem}
By \eqref{6} and \eqref{18}, we note that
\begin{equation}
\beta_{n,\lambda}^{(s)}(x)\sim\bigg(\frac{(e_{\lambda}(t)-1)^{s}}{t^{s}},t\bigg)_{\lambda} \label{28}
\end{equation}
From \eqref{20}, \eqref{23} and \eqref{28}, we have
\begin{equation}
B_{n,\lambda}^{(k)}(x)=\sum_{m=0}^{n}C_{n,m}\beta_{m,\lambda}^{(s)}(x),\label{29}
\end{equation}
where
\begin{align}
C_{n,m}\ &=\ \frac{1}{m!}\Bigg\langle \frac{\mathrm{Li}_{k,\lambda}(1-e_{\lambda}(-t))}{e_{\lambda}(t)-1}\frac{(e_{\lambda}(t)-1)^{s}}{t^{s}}t^{m}\ \Bigg|\ (x)_{n,\lambda}\Bigg\rangle_{\lambda} \label{30} \\
&= \binom{n}{m} \Bigg\langle \frac{\mathrm{Li}_{k,\lambda}(1-e_{\lambda}(-t))}{e_{\lambda}(t)-1}\frac{(e_{\lambda}(t)-1)^{s}}{t^{s}} \Bigg|\ (x)_{n-m,\lambda}\Bigg\rangle_{\lambda}    \nonumber \\
&=\ \binom{n}{m}\sum_{l=0}^{n-m}\frac{S_{2,\lambda}(l+s,s)}{\binom{l+s}{s}l!} \Bigg\langle \frac{\mathrm{Li}_{k,\lambda}(1-e_{\lambda}(-t))}{e_{\lambda}(t)-1}t^{l}\ \Bigg|\ (x)_{n-m,\lambda}\Bigg\rangle_{\lambda}\nonumber\\
&=\ \binom{n}{m}\sum_{l=0}^{n-m}\frac{S_{2,\lambda}(l+s,s)}{\binom{l+s}{s}}\binom{n-m}{l} \Bigg\langle \frac{\mathrm{Li}_{k,\lambda}(1-e_{\lambda}(-t))}{e_{\lambda}(t)-1}\ \Bigg|\ (x)_{n-m-l,\lambda}\Bigg\rangle_{\lambda} \nonumber \\
&=\ \binom{n}{m}\sum_{l=0}^{n-m}\frac{S_{2,\lambda}(l+s,s)}{\binom{l+s}{s}}\binom{n-m}{l}B_{n-m-l,\lambda}^{(k)}.\nonumber
\end{align}
Therefore, by \eqref{29} and \eqref{30}, we obtain the following theorem.
\begin{theorem}
For $n\ge 0$, we have
\begin{displaymath}
B_{n,\lambda}^{(k)}(x)=\sum_{m=0}^{n}\binom{n}{m}\bigg(\sum_{l=0}^{n-m}\frac{\binom{n-m}{l}}{\binom{l+s}{s}}S_{2,\lambda}(l+s,s)B_{n-m-l,\lambda}^{(k)}\bigg)\beta_{m,\lambda}^{(s)}(x).
\end{displaymath}
\end{theorem}
For $n\ge 0$, we let
\begin{displaymath}
\mathbb{P}_{n}=\big\{p(x)\in\mathbb{C}[x]\ \big|\ \deg p(x)\le n\big\}.
\end{displaymath}
Then $\mathbb{P}_{n}$ is an $(n+1)$-dimensional vector space over $\mathbb{C}$. \par
From \eqref{11}, we note that $d_{n,\lambda}(x)\sim ((1-t)e_{\lambda}(t),t)_{\lambda},\ (n\ge 0)$.
For $p(x)\in\mathbb{P}_{n}$, we let
\begin{equation}
p(x)=\sum_{l=0}^{n}C_{l}d_{l,\lambda}(x). \label{31}
\end{equation}
By \eqref{17}, we have
\begin{align}
\langle (1-t)e_{\lambda}(t)t^{m}\ |\ p(x)\rangle_{\lambda}\ &=\ \sum_{l=0}^{n}C_{l}\big\langle (1-t)e_{\lambda}(t)t^{m}\ \big|\ d_{l,\lambda}(x)\rangle_{\lambda}\label{32}\\
&=\ \sum_{l=0}^{n}C_{l}m!\delta_{m,l}\ =\ C_{m}m!,\nonumber
\end{align}
where $0\le m\le n$. \par
Therefore, by \eqref{31} and \eqref{32}, we obtain the following theorem.
\begin{theorem}
For $p(x)\in\mathbb{P}_{n}$, we have
\begin{displaymath}
p(x)=\sum_{l=0}^{n}C_{l}d_{l,\lambda}(x),
\end{displaymath}
where
\begin{displaymath}
C_{l}=\frac{1}{l!}\big\langle (1-t)e_{\lambda}t^{l}\big|p(x)\big\rangle_{\lambda}.
\end{displaymath}
\end{theorem}
Let $p(x)=B_{n,\lambda}^{(k)}(x)\in\mathbb{P}_{n}$. Then we have
\begin{equation}
B_{n,\lambda}^{(k)}(x)=\sum_{l=0}^{n}C_{l}d_{l,\lambda}(x),\label{33}
\end{equation}
where
\begin{align}
C_{l}\ &=\ \frac{1}{l!}\big\langle (1-t)e_{\lambda}(t)t^{l}|B_{n,\lambda}^{(k)}(x)\big\rangle_{\lambda}\label{34} \\
&=\ \binom{n}{l}\big\langle (1-t)e_{\lambda}(t)\big|B_{n-l,\lambda}^{(k)}(x)\big\rangle_{\lambda}\nonumber \\
&=\ \binom{n}{l}\big\langle (1-t)|B_{n-l,\lambda}^{(k)}(x+1)\big\rangle_{\lambda} \nonumber\\
&=\ \binom{n}{l}\big\langle 1\big|B_{n-l,\lambda}^{(k)}(x+1)\big\rangle_{\lambda}-\binom{n}{l}(n-l)\big\langle 1|B_{n-l-1,\lambda}^{(k)}(x+1)\big\rangle_{\lambda}\nonumber\\
&=\ \binom{n}{l}B_{n-l,\lambda}^{(k)}-n\binom{n-1}{l}B_{n-l-1,\lambda}^{(k)}(1). \nonumber
\end{align}
Thus, by \eqref{33} and \eqref{34}, we get
\begin{equation}
B_{n,\lambda}^{(k)}(x)=\sum_{l=0}^{n}\bigg(\binom{n}{l}B_{n-l,\lambda}^{(k)}-n\binom{n-1}{l}B_{n-l-1,\lambda}^{(k)}\bigg)d_{l,\lambda}(x).\label{35}
\end{equation}
From \eqref{11}, we note that $d_{n,r}^{(r)}(x)\sim\big((1-t)^{r}e_{\lambda}(t),t\big)_{\lambda}$. \par
Let us assume that
\begin{equation}
p(x)=\sum_{m=0}^{n}C_{m}^{(r)}d_{m,\lambda}^{(r)}(x)\in\mathbb{P}_{n}. \label{36}
\end{equation}
Then, by \eqref{17}, we get
\begin{align}
\big\langle (1-t)^{r}e_{\lambda}(t)t^{m}\big|p(x)\big\rangle_{\lambda}\ &=\ \sum_{l=0}^{n}C_{l}^{(r)}\big\langle (1-t)^{r}e_{\lambda}(t)t^{m}\ \big|\ d_{l,\lambda}^{(r)}(x)\big\rangle_{\lambda}\label{37} \\
&=\ m!C_{m}^{(r)},\quad (0\le m\le n).\nonumber
\end{align}
Therefore, by \eqref{36} and \eqref{37}, we obtain the following theorem.
\begin{theorem}
    For $n\ge 0$, we have
    \begin{displaymath}
        p(x)=\sum_{m=0}^{n}C_{m}^{(r)}d_{m,\lambda}^{(r)}(x),
    \end{displaymath}
    where
    \begin{displaymath}
        C_{m}^{(r)}=\frac{1}{m!}\big\langle (1-t)^{r}e_{\lambda}(t)t^{m}\ \big|\ p(x)\big\rangle.
    \end{displaymath}
\end{theorem}
We let
\begin{equation}
d_{n,\lambda}(x)=\sum_{m=0}^{n}C_{m}^{(r)}d_{m,\lambda}^{(r)}(x),   \label{38}
\end{equation}
where
\begin{align}
C_{m}^{(r)}\ &=\ \frac{1}{m!}\big\langle (1-t)^{r}e_{\lambda}(t)t^{m}\ \big|\ d_{n,\lambda}(x)\big\rangle_{\lambda}\label{39} \\
&=\ \binom{n}{m}\big\langle (1-t)^{r}\ \big|\ d_{n-m,\lambda}(x+1)\big\rangle_{\lambda}\nonumber \\
&=\ \binom{n}{m}\sum_{j=0}^{r}\binom{r}{j}(-1)^{j}(n-m)_{j}\big\langle 1\big|d_{n-m-j,\lambda}(x+1)\rangle_{\lambda}\nonumber \\
&=\ \binom{n}{m}\sum_{j=0}^{r}\binom{r}{j}\binom{n-m}{j}(-1)^{j}d_{n-m-j,\lambda}(1)j!. \nonumber
\end{align}
By \eqref{38} and \eqref{39}, we get
\begin{equation}
d_{n,\lambda}(x)=\sum_{m=0}^{n}\binom{n}{m}\bigg(\sum_{j=0}^{r}\binom{r}{j}\binom{n-m}{j}j!(-1)^{j}d_{n-m-j,\lambda}(1)\bigg)d_{m,\lambda}^{(r)}(x),\quad(n\ge 0). \label{40}
\end{equation}
Let us take $p(x)=B_{n,\lambda}^{(k)}(x)\in\mathbb{P}_{n}$, $(n\ge 0)$. Then, by Theorem 7, we get
\begin{equation}
B_{n,\lambda}^{(k)}(x)=\sum_{m=0}^{n}C_{m}^{(r)}d_{m,\lambda}^{(r)}(x),\quad(n\ge 0),\label{41}
\end{equation}
where
\begin{align}
C_{m}^{(r)}\ &=\ \frac{1}{m!}\big\langle (1-t)^{r}e_{\lambda}(t)t^{m}\big|B_{n,\lambda}^{(k)}(x)\big\rangle_{\lambda} \label{42} \\
&=\ \binom{n}{m}\big\langle (1-t)^{r}\big|B_{n-m,\lambda}^{(k)}(x+1)\big\rangle_{\lambda}\nonumber \\
&=\ \binom{n}{m}\sum_{j=0}^{r}\binom{r}{j}(-1)^{j}\binom{n-m}{j}j!  \big\langle 1\big|B_{n-m-j}^{(k)}(x+1)\big\rangle_{\lambda}\nonumber \\
&=\ \binom{n}{m}\sum_{j=0}^{r}\binom{r}{j}(-1)^{j}\binom{n-m}{j}j!B_{n-m-j,\lambda}^{(k)}(1). \nonumber
\end{align}
Therefore, by \eqref{41} and \eqref{42}, we obtain the following theorem.
\begin{theorem}
    For $n\ge 0,\ r\in\mathbb{N}$, we have
    \begin{displaymath}
        B_{n,\lambda}^{k}(x)=\sum_{m=0}^{n}\binom{n}{m}\bigg(\sum_{j=0}^{r}\binom{r}{j}(-1)^{j}\binom{n-m}{j}j!B_{n-m-j,\lambda}^{(k)}(1)\bigg)d_{m,\lambda}^{(r)}(x).
    \end{displaymath}
\end{theorem}

\section{Conclusion}
The study of degenerate versions of some special polynomials and numbers, which was initiated by Carlitz, regained interests of some mathematicians and many interesting results were found not only in their arithmetical and combinatorial aspects but also in applications to differential equations, identities of symmetry and probability theory. \par
Recently, the `$\lambda$-umbral calculus' was developed by the motivation that what if the usual exponential function is replaced by the degenerate exponential functions in the generating function of a Sheffer sequence. This question led us to the introduction of the concepts like $\lambda$-linear functionals,  $\lambda$-differential operators and $\lambda$-Sheffer sequences. \par
In this paper, we studied the degenerate poly-Bernoulli polynomials, which is a $\lambda$-Sheffer sequence and hence a degenerate Sheffer sequence, by using three different formulas, namely a formula about representing a $\lambda$-Sheffer sequence by another, a formula coming from the generating functions of  $\lambda$-Sheffer sequences and a formula arising from the definitions for $\lambda$-Sheffer sequences.
Then, among other things, we represented the degenerate poly-Bernoulli polynomials by higher-order degenerate Bernoulli polynomials and by higher-order degenerate derangement polynomials. \par
It is one of our future projects to continue to investigate the degenerate special numbers and polynomials by using the recently developed $\lambda$-umbral calculus.

\bigskip

\noindent{\bf{Acknowledgments}}\\
The authors would like to thank Jangjeon Research Institute for Mathematical Sciences for the support of this research.

\bigskip

\noindent{\bf{Funding}}\\
Not applicable.

\bigskip

\noindent{\bf{Availability of data andmaterials}}\\
Not applicable.

\bigskip

\noindent{\bf{Ethics approval and consent to participate}}\\
All authors reveal that there is no ethical problem in the production of this paper.

\bigskip

\noindent{\bf{Competing interests}}\\
The authors declare that they have no competing interests.

\bigskip

\noindent{\bf{Consent for publication}}\\
All authors want to publish this paper in this journal.

\bigskip

\noindent{\bf{Authors' contributions}}\\
TK and DSK conceived of the framework and structured the whole paper; TK and DSK wrote the paper; JK and HL checked the results of the paper and typed the paper; DSK and TK completed the revision of the article. All authors have read and agreed to the published version of the manuscript.

\bigskip

\noindent{\bf{Author details}}\\


\begin{thebibliography}{9}
\bibitem{1}
Carlitz, L. \emph{Degenerate Stirling, Bernoulli and Eulerian numbers,} Utilitas Math. \textbf{15} (1979), 51–-88
\bibitem{2}
Dere, R.; Simsek, Y. \emph{Applications of umbral algebra to some special polynomials,}
Adv. Stud. Contemp. Math. (Kyungshang) \textbf{22} (2012), no. 3, 433--438.
\bibitem{3}
Gzyl, H. \emph{Canonical transformations, umbral calculus, and orthogonal theory,} J. Math. Anal. Appl.  \textbf{111} (1985), no. 2, 547–-558.
\bibitem{4}
Kholodov, A. N. \emph{The umbral calculus and orthogonal polynomials,} Acta Appl. Math. \textbf{19} (1990), no. 1, 1–-54
\bibitem{5}
Kim, D. S.; Kim, T. \emph{A note on a new type of degenerate Bernoulli numbers,} Russ. J. Math. Phys. \textbf{27} (2020), no. 2, 227–-235.
\bibitem{6}
Kim, D. S.; Kim, T. \emph{Sheffer sequences and $\lambda$-Sheffer sequences,} J. Math. Anal. Appl. \textbf{493} (2021), no. 1, 124521, 21 pp.
\bibitem{9}
Kim, D. S.; Kim, T. \emph{A note on polyexponential and unipoly functions,} Russ. J. Math. Phys. \textbf{26}  (2019), no. 1, 40–-49.
\bibitem{8}
Kim, T.; Kim, D. S. \emph{Note on the degenerate gamma Function,} Russ. J. Math. Phys. \textbf{27}  (2020), no. 3, 352–-358.
\bibitem{11}
Kim, T.; Kim, D. S. \emph{Degenerate polyexponential functions and degenerate Bell polynomials,} J. Math. Anal. Appl. \textbf{487} (2020), no. 2, 124017, 15 pp.
\bibitem{12}
Kim, T.; Kim, D. S. \emph{Some identities on derangement and degenerate derangement polynomials,}  Advances in mathematical inequalities and applications, 265–-277, Trends Math., Birkh\"auser/Springer, Singapore, 2018.
\bibitem{13}
Kim, T.; Kim, D. S.; Dolgy, D. V.; Kwon, J. \emph{Some identities of derangement numbers,} Proc. Jangjeon Math. Soc.  \textbf{21} (2018), no. 1, 125–-141.
\bibitem{15}
Kim, T.; Kim, D. S.; Jang, G.-W.; Kwon, J. \emph{A note on some identities of derangement polynomials,} J. Inequal. Appl. 2018, Paper No. 40, 17 pp.
\bibitem{7}
Kim, T.; Kim, D. S.; Kim, H.-Y.; Lee, H.; Jang, L.-C. \emph{Degenerate poly-Bernoulli polynomials arising from degenerate polylogarithm,} Adv. Difference Equ. 2020, Paper No. 444, 9 pp
\bibitem{14}
Kim, T.; Kim, D. S.; Kwon, H.-I. Jang, L.-C. \emph{Fourier series of sums of products of $r$-derangement functions,} J. Nonlinear Sci. Appl. \textbf{11} (2018), no. 4, 575–-590.
\bibitem{10}
Kim, T.; Kim, D. S.; Lee, H.; Jang, L.-C. \emph{A note on degenerate derangement polynomials and numbers,} arXiv:2011.08535. https://arxiv.org/abs/2011.08535
\bibitem{16}
Lee, D. S.; Kim, H. K. \emph{On the new type of degenerate poly-Genocchi numbers and polynomials,} Adv. Difference Equ. 2020, Paper No. 431, 15 pp.
\bibitem{17}
Ray, N. \emph{Extensions of umbral calculus: penumbral coalgebras and generalised Bernoulli numbers,} Adv. in Math. \textbf{61} (1986), no. 1, 49–-100.
\bibitem{18}
Roman, S. \emph{The umbral calculus,} Pure and Applied Mathematics, 111. Academic Press, Inc. [Harcourt Brace Jovanovich, Publishers], New York, 1984.
\bibitem{19}
Rota, G.-C.; Taylor, B. D. \emph{The classical umbral calculus,} SIAM J. Math. Anal. \textbf{25} (1994), no. 2, 694–-711.
\bibitem{20}
Ueno, K. \emph{General power umbral calculus in several variables,} J. Pure Appl. Algebra \textbf{59} (1989), no. 3, 299–308.
\bibitem{21}
Wilson, B. G.; Rogers, F. J. \emph{Umbral calculus and the theory of multispecies nonideal gases,} Phys. A \textbf{139} (1986), no. 2-3, 359–-386.
\end{thebibliography}
\end{document}